# Entropy And Vision

Rami Kanhouche[1]

26-6-2006

*Abstract*

In vector quantization the number of vectors used to construct the codebook is always an undefined problem, there is always a compromise between the number of vectors and the quantity of information lost during the compression. In this text we present a minimum of Entropy principle that gives solution to this compromise and represents an Entropy point of view of signal compression in general. Also we present a new adaptive Object Quantization technique that is the same for the compression and the perception.

## I. Introduction

Vector quantization [2] is the process of classification of a group of vectors $\{v_i\}_{i:=0...M-1}, v_i \in R^k$ into a smaller number of vectors group $\{v'_i\}_{i:=0,N-1}, v'_i \in R^k$, $N \ll M$, the classification of each vector $v_i$ is the mapping $I(i):\{0,1,2,...M-1\} \to \{0,1,2,...N-1\}$. By this, a vector $v_i$ would be classified by the vector $v'_{I(i)}$. The mapping is always from the vector $v_i$ to the « nearest » vector $v'_j$, so that

$$I(i) = \arg\min_j \|v_i - v'_j\|. \qquad (1)$$

The classification vectors or what is called the codebook should minimize the distance between each vector and its class vector, and for that the group of vectors $\{v'_i\}$ should take values minimizing

$$\text{Min} \quad D(\{v_i\},\{v'_j\}) := \sum_{i=0}^{M-1} \|v_i - v'_{I(i)}\|. \qquad (2)$$
$$\text{Subject to } \{v'_i\}$$

The process of vector quantization is used in many applications. In this text we are interested of its informational aspect. We point by the informational aspect its capability to resume or « compress » a group of vectors and provide a reasonable quantification of this group distribution. This ability to "compress" –or summarize-

---

[1] PhD Student at Lab. CMLA, École Normale Supérieure de Cachan, 61, avenue du Président Wilson, 94235 CACHAN Cedex, France. Phone: +33-1-40112688, mobile: 33-6-62298219, fax: +33-1-47405901.E-mail : rami.kanhouche@cmla.ens-cachan.fr, kanram@free.fr.
Web: http://www.cmla.ens-cachan.fr/Utilisateurs/kanhouch

R. Kanhouche

which got a mental-or qualitative importance in application when a person needs to have an intelligent mathematical view of the data; got also a quantitative importance that find a direct application in image compression, and data compression in general. In what concerns image compression, usually the image is divided into $M$ blocks of pixels, each corresponding to a vector $v_i$. We found that the percentage in which the data size reduction happens is

$$T(\{v_i\},\{v'_j\}) := \frac{Pk}{\log_2 N + \frac{PkN}{M}} \quad (3)^2$$

where $k$ is the dimension of the vector space, and $P$ is the number of bits needed to store each value $v_i^h$ composing the vectors $\{v_i\}$ $v_i := [v_i^h]_{h:=0...k-1}$.

In real application the numbers $M$, $P$ and $k$ are usually predefined. On the other hand the number $N$ imply a direct decision over the values of each of $\{v'_i\}$, and $D(\{v_i\},\{v'_j\})$. The selection of $N$, is always a compromise between the efficiency of the data reduction $(T)$, and the sincerity of the classification $\left(\frac{1}{D}\right)$. In this text we will answer the question about the appropriate selection of $N$ as a solution to this compromise.

## II. Quantization From Entropy perspective

It is very natural to say that quantization tries to minimize the Entropy [1] of the signal in an effort to be able have an efficient data size reduction. To formalize this we can transfer the problem in $(3)$ as

$$\text{Min } \delta.a - b.\sum_i P(v'_i).\log(P(v'_i)) \quad (4)$$

$$\text{Subject to } \{v'_i\}, \delta, N \leq M : \forall v_i, \exists v'_j : \|v_i - v'_j\| < \delta$$

Where the probability of $v'_i$ is the number of vectors $v_i$ classified as $v'_i$,

$$P(v'_i) := \frac{|\{v_j : I(j) = i\}|}{N}.$$

What is advantageous in this formalization is that it is an accurate expression of the compromise as we come to explain in the introduction. The difference between (4) and (3) is that the number $N$ is not fixed. The values $a$ and $b$ are constants that balance between data size reduction, and the accuracy of the solution, or in other words, they balance between $T$ and $D$. To help interpreting the formalization in (4) one could observe that a very accurate matching would give a minimal $\delta$, with a maximum of Entropy. On the other side, for a given $N$ there is an optimum minimum of Entropy with a given $\delta$ that could make the solution unrepresentative.

---

[2] The explanation of this relation is that the percentage of compression is the number of bits needed to store each image vector (block) before compression divided on the needed number of bits needed to store that block after. The later is the number of bits needed to store the class number of that block $(\log_2 N)$, and with it the actual number of bits needed to store the code book itself $(PkN)$ distributed over all blocks $(M)$.

This contradiction is only technical, meaning, that when presenting a given number $N$ the same process minimizing the distance according to (3) is actually minimizing the Entropy. And from that it is quite reasonable to wonder if there is another form that is more expressive of the quantization process.

This goes very profoundly in the manner in which the visual system estimates the data, which as we will argue in the next section is an informational manner, which put in synchronization the quality and the quantity.

## III. The Submerging Entropy Estimate

The Submerging Entropy Estimate is defined as the following

$$SEE(\{v_i\}) = \begin{cases} \min_{\text{subject to } \{v'_i\}} \sum_i P(v'_i)\left[\log\left(\frac{1}{P(v'_i)}\right) + 1 + \frac{1}{P(v'_i)}\frac{|\Delta(v'_i)|}{|\{v_i\}|}SEE(\Delta(v'_i))\right] & \text{if } |\{v_i\}| > 1 \\ 0 & \text{if } |\{v_i\}| \leq 1 \end{cases}$$

(5)

where $\Delta(v'_i) := \{v_j - v'_i : v_j - v'_i \neq 0, I(j) = i\}$

The next theorem is a direct explanation of the previous definition

**Definition 1** We call an $N$-level auto generative events group, any group of events for which each event do generate a new events group, in a recursive manner for $N$ levels. We note $G_t^i$ for the events group generated at the level $t$ from the event $i$ at the level $t-1, 1 \leq t < N$, with the source events group noted as $G_0$.

**Theorem 1** The entropy of $N$-level auto generative events group is equal to

$$H_T^N(G_0)$$

where

$$H_T^N(G_t) := \begin{cases} \sum_{i \in G_t} P(i)\left[\log\left(\frac{1}{P(i)}\right) + H_T^N(G_{t+1}^i)\right] & \text{if } t < N-1 \\ H(G_t) & \text{if } t = N-1 \end{cases}$$

and $H(G) := \sum_{i \in G} P(i).\log\left(\frac{1}{P(i)}\right)$.

The previous structure represents what we can call a non-stationary Active Markov Chain. Figure 1 sketches the topology of this kind of chain compared to the other types of Markov chains. Active Markov Chains are the generalization of static Markov Chains; here are some of the properties of the process following an Active Markov Chains model:

- It got entropy.

- It is Active, meaning it got an adaptive model minimizing its entropy.

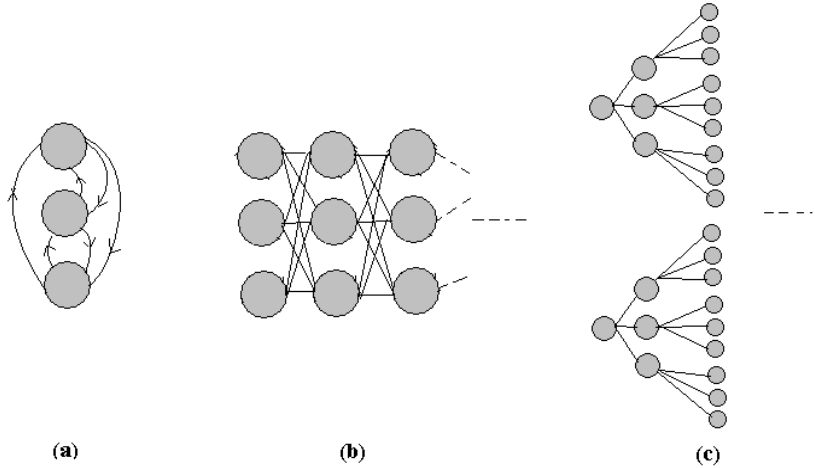

**Figure 1, (a) stationary, (b) non-stationary static, (c) non-stationary active.**

The definition (5) has a recursive nature. Actually it defines a hierarchical Entropy Tree that adapt to the data in a way that realize an optimum compression. It depends as its structure explains on the principle of **"nothing is left behind"**. Meaning, that a global view for the multilevel compression of the vector space $\{v_j\}$ is always considered, even that for real application only part of the top-level nodes –down to a given depth- would be preserved to reconstruct the signal. Each node got a quantity value expressing its data size, and by that when selecting the nodes that will participate in the construction we have always a criteria to differentiate between the perceptional quality of each configuration. Also the advantage in here is that *N* is taken out of the equation, and instead, we got a global point of view of all possible combinations of Vector Quantization which are independent from *N,* and at the same time depending on a consistent value which is the Entropy.

From Theorem 1 we can see the logical structure of (5), the difference in (5) is that not every event would produce events at the next level. More precisely, when there is a zero mach then the population of the next level is altered (reduced), this is taken into consideration by the correcting factor $\frac{1}{P(v'_i)} \frac{|\Delta(v'_i)|}{|\{v_i\}|}$. Also 1 is added to count for coding the situation when an event stops at a level.

The convergence of (5) according to **Theorem 1** is guided by two contradictions:
1- The need to stop the recursion by approaching the elements $\{v_i\}$ when there is no entropy left.
2- The need to extend the recursion to have a good compression at each step when there is a possibility to that.

This multilevel view of the data is not new, it was implicitly considered for Vector Quantization by [3], [6]. What is new in the definition of SEE is the explicit expression of this point of view in the Entropy sense. In other words it is a quantization of the vector quantization it self. It also could help providing a unique gate to all kind of possible quantization scenario, not only vectors of square parts- blocks- of the image.

## IV. Toward a global framework for the compression

In the previous section we introduced an Entropy point of view of the compression using Vector Quantization. At the same time there is many different methods for image compression, each with its advantages and disadvantages. Fractal image compression is one of the most intriguing image compression methods. In the original fractal image compression, which could be considered is a "more intelligent" compression; only the relations between different parts of the image are saved to construct the entire image without the need of any of the actual data present in each part. The origin of the fractal compression idea is the definition of a class of operators called contractive operators having the property:

$$\forall x_0 \in R^N, \exists x : \lim f(...f(x_0)) = x$$

This process, of fractal compression, is to build the operator $f$ so that it is contractive toward $x$ -the image. There is many ways to define this operator in a function of the image subject of the compression. One of these ways is to divide the image recursively into overlapping triangles –or squares, and after, to project each triangle- or square- on a triangle-or square- of smaller size [11].

Also as in the Vector Quantization compression this could have multiple solutions. In other words, there could be multiples of $f$, each contractive toward the same $x$. By the same way, as in Vector Quantization, this ambiguity could find in the Entropy expression a reasonable quantifying measure.

For an example, and in the vector quantization sense, the vector that we are trying to get quantized contains the difference between the projected part and the projected-on part, the index number of the projected-on part –i.e. triangle number, when using triangles- and the other parameters like the rotational parameters and the shrinking parameter.

Of course the relation in here is not direct as in the Vector Quantization situations. Especially, in what connect the quantity of the information to the visual aspect of the reconstruction -decoding. From which we can understand that the Entropy expression of the vector quantization in (5) is not sufficient as a general model. The development of such a model starts from observing that all kind of compression methods follows the following primitive steps:

1- Compression
2- Compressed situation Data.
3- Decompression

This scheme, as simple as it is, is very fundamental, and when expressed in algebra it becomes as

$$x = f(x). \qquad (6)$$

In matrix algebra this could be expressed as

$$F x = x, \qquad (7)$$

where $F$ is the compression matrix of size $N \times N$, $N$ is the signal size.
The global model that we will call the Compression Gate take the following form:

$$\min Entropy(f, x) \qquad (8)$$
$$\text{subject to } x = f(x)$$

The Entropy term $(f,x)$ is by definition the quantity of information needed to obtain $f$, and $x$. The previous model is not a magical tool that realizes the optimum compression; it is just a chance to realize that. What is advantageous according to the (8) relation is that $f$ is not only a theoretical representation of other methods for the compression, but instead it could be expressed according to (7) to represent the calculation of many compression methods.

For an example in the compression by Fourier Transform it take the form:
$$f(x) := \mathbf{e}^H X(w_1,w_2) = \mathbf{e}^H \mathbf{e} x$$

where $\mathbf{e}$ is the Fourier Transformation matrix and $X(w_1,w_2)$ is the image in the frequency domain. The Entropy of $X(w_1,w_2)$, which is in this context equivalent to $f$, has the maximum quantity of information in the low frequencies components. In other words when considering a vector quantization of $X(w_1,w_2)$ into blocks $\{X_i(w_1,w_2)\}_i$:

$$SEE(\{\text{Re}(X_i(w_1,w_2))\}) + SEE(\{\text{Im}(X_i(w_1,w_2))\}),$$

the end nodes with the lower quantity of information correspond to those dealing with the high frequencies. One could roughly say that the Fourier transformation is an implicit search into the *SEE* operator.

For the compression using the wavelet, which is superior to the Fourier compression in compression efficiency, the same kind of discussion does apply.

## V. Object and Entropy

Until here we have considered quantization, and eventually the compression, using only vectors of equal size. Segmenting the image into blocks of equal size produces these vectors. A more adaptive form for modeling the image would be not to limit our self to this constraint. In [12] we proposed an Object point of view of the image. In this context the vector of values become a surface of points into a binary space. To expose the problem in this form it is necessary to be able to define an adaptive configuration of the quantization "vector"-Object.

The Object point of view as it was presented in [12] was presented in the context of Object Perception context. On the other hand there is an inherent relation between the compression and the perception that we had reported in the past, and was reported also in [10], and we are clarifying more in this text. In fact the relation (8) is a clear proof of that relation, since we have showed how object perception and compression performs in the same space; the space of $x = f(x)$.

To start we will consider the image $S(x,y):[0,1,\ldots N_0-1]\times[0,1,\ldots N_1-1]\to R^+$, and more generally any signal $S(x'):\prod_{k=0}^{d-1}[0,1,\ldots,N_k-1]\to R$, where $d$ is the signal dimension, and $\{N_k\}_{k:=0,1\ldots d-1}$ is the signal size in each dimension, as an Object of $\prod_{k=0}^{d-1} N_k$ points into a binary space. This is formalized according to

$$I_S(x): \prod_{k=0}^{d}[0,1,\ldots,N_k-1] \to \{0,1\}.$$

$$I_S(x) := \begin{cases} 1 & if \quad S([x_0,x_1,\ldots,x_{d-1}]) = x_d \\ 0 & if \quad S([x_0,x_1,\ldots,x_{d-1}]) \neq x_d \end{cases} \qquad (9)$$

with $N_d$ is the maximum discreet value the signal could have.

In [12] the Objects was presumed of infinite support. This is not appropriate for real application. For that the distance between two Objects with the same number of points $N$, would be defined as a function of two values as the following

$$\chi_{\alpha,\gamma}(I_1,I_2) := \sum_{x \in \mathfrak{I}_1, y \in \mathfrak{I}_1} \left| dist(x,y) - \alpha \, dist(\Omega_\gamma(x), \Omega_\gamma(y)) \right| \qquad (10)$$

$\Omega_\gamma : \mathfrak{I}_1 \to \mathfrak{I}_2$, where $\mathfrak{I}_u := \{x : I_u(x) \neq 0\}$

The value $\alpha$ is the projection parameter, and the $\gamma$ value is the combination parameter. In [12] we have showed how the combinatorial space of $\Omega : \mathfrak{I}_1 \to \mathfrak{I}_2$, could be projected on the line. And for that we noted in here $\Omega_\gamma$ as the permutation corresponding to the value $\gamma$ on the combinatorial line. In matrix notation we would have

$$\chi_{\alpha,\gamma}(I_1,I_2) = \sum_{ij} \left[ D_2 - \frac{1}{\alpha} W_\gamma^T D_1 W_\gamma \right]_{ij} \qquad (11)$$

where $D_u$ is the auto-distance matrix corresponding to the signal $I_u$, defined as

$$D_u := [d_{ij}]_{ij:=0,\ldots,N-1}, \; d_{ij} := dist(x_i, x_j). \qquad (12)$$

And $W_\gamma$ is the permutation matrix defined as

$$W_\gamma := [w_{ij}]_{ij:=0,\ldots N-1}, w_{ij} := \begin{cases} 1 & if \quad j = \Omega(i) \\ 0 & if \quad j \neq \Omega(i) \end{cases}. \qquad (13)$$

Before that we get into the Entropy expression in the Object sense we will need few definitions.

**Definition 2** an Object is a finite group of points $O := \{x_i \in R^{d+1} : i \in \{0,\ldots,N-1\}\}$, we call $C_O := \frac{1}{|O|} \sum_i x_i, x_i \in O$ as the center of $O$.

**Definition 3** the distance between each two Objects of the same number of points is defined as $\chi_{\alpha,\gamma}(O_1,O_2) := \sum_{x \in O_1, y \in O_1} \left| dist(x,y) - \alpha \, dist(\Omega_\gamma(x), \Omega_\gamma(y)) \right|$
$\Omega_\gamma : O_1 \to O_2$.

**Definition 4** we define the difference between two Objects $O_1, O_2$ as
$diff_{\alpha,\gamma}(O_1,O_2) := \left\{ x_i : x_i = dist\left(y_i - C_{O_1}, \alpha.(\Omega_\gamma(y_i) - C_{O_2})\right), y_i \in O_1 \right\}$.

**Definition 5** we call a Generalization Operator $G_N$ any Object operator that produce from any given Object, and for any $N$ a new Object having an $N$ number of points, so that $\forall O, N_2 : |G_{N_2}(O)| = N_2$.

**Definition 6** we call the space $P(O)$ any space of objects belonging to $O$ so that $\bigcup P(O) = O$, $\forall O_1 \in P(O), O_2 \in P(O), O_1 \cap O_2 = \phi$.

The Entropy expression in the Object sense would become

$$SEE(\{O_i\}) = \begin{cases} \min_{\text{subject to } \{O_i'\}} \sum_i P(O_i') \cdot \left[ \log\left(\frac{1}{P(O_i')}\right) + 1 + \frac{1}{P(O_i')} \frac{|\Delta(O_i')|}{|\{O_i\}|} SEE(\Delta(O_i')) \right] & \text{if } |\{O_i\}| > 1 \\ 0 & \text{if } |\{O_i\}| \leq 1 \end{cases}$$

(14)

where

$$\Delta(O_i') := \left\{ \text{diff}_{\alpha,\gamma}\left(G_{|O_j|}(O_i'), O_j\right) : \text{diff}_{\alpha,\gamma}\left(G_{|O_j|}(O_i'), O_j\right) \neq \mathbf{0}, (\alpha, \gamma, i) = \arg\min\left(\chi\left(G_{|O_j|}(O_i'), O_j\right)\right) \right\}$$

In the previous formalization, and for simplicity, we did not take the elements $\alpha, \gamma$ into consideration. Meanwhile, this could be done by increasing the dimension of each Object in $\Delta(O_i')$. For a given Object –or image, the problem $SEE(\{O_i\})$ is also subject to the problem

$$\min_{\text{subject to } P(O)} SEE(P(O)) \qquad (15)$$

## VI. Conclusion

We had the chance to provide an Entropy point of view to both of quantization and perception. In what concern Vector Quantization we provided a more logical manner for the decision on the number of the quantization vectors. The principle of "nothing is left behind" is very fundamental. It represents a new technique to provide an in depth understanding of many compression techniques, and helps having a unified perspective for both of the compression without lost with the compression with lost.

The new concept of Active Markov Chains is very fundamental, from one side it could be used in any compression (identification) scenario, and at the same time it defines stochastic processes that are of a great degree of adaptability.

Also we have presented an Entropy expression for Object Analysis, which could provide an enhanced compression ratio. Of course the compression according to (15), or (14) is not an optimal compression, but explain how there is always a way to minimize the Entropy of:

$$x = f(x).$$